\documentclass[a4paper,10pt]{article}
\usepackage{color}
\usepackage{cite}
\usepackage{amsmath,amsthm,amssymb,amsfonts}
\usepackage{eucal}
\usepackage[all]{xy}
\usepackage{color}
\usepackage{fancyhdr}
\usepackage{lastpage}
\usepackage{graphicx}
\usepackage{bm,bbm}
\usepackage{ascmac,paralist}
\usepackage{url}
\usepackage{array,arydshln}
\usepackage{tikz}
\usetikzlibrary{intersections,calc,arrows.meta}

\newtheorem{theorem}{Theorem}[section]

\newtheorem{corollary}[theorem]{Corollary}
\newtheorem{lemma}[theorem]{Lemma}
\newtheorem{proposition}[theorem]{Proposition}

\newtheorem{assumption}[theorem]{Assumption}

\title{A Spectral Analysis of \\
The Correlated Random Walk}
\author{Yusuke Ide\thanks{Department of Mathematics, College of Humanities and Sciences, Nihon University, Setagaya, Tokyo, 156-8550, Japan, E-mail: ide.yusuke@nihon-u.ac.jp},
Akihiro Narimatsu\thanks{Center for Mathematical and Data Scienses, The University of Fukuchiyama, Fukuchiyama, Kyoto, 620-0886, Japan, E-mail: narimatsu-akihiro@fukuchiyama.ac.jp}
}
\date{}
\begin{document}
\maketitle
\begin{abstract}
In this paper, we consider a spectral analysis of the Correlated Random Walk (CRW) on the path. We apply an analytical method for the Quantum Walk to CRW. For the isospectral coin cases, we obtain all of the eigenvalues and the corresponding eigenvectors of the time evolution operator of CRW, and also obtain the limiting distribution. 
\end{abstract}
\section{Introduction}

The classical Random Walk (RW) is one of the most important models in the study of probability theory and related fields. Furthermore, a more general model called the Correlated Random Walk (CRW) or the Persistent Random Walk has been studied by various authors.  \cite{alla,bohm,bohm2,persis,rwinbio,insect,kiumi,kks1,konno2,someexpli}. 

The motion of CRW depends on the movement of the previous step. This feature is regarded as the random walker to have multiple states at each vertex. The similar feature is seen in the Quantum Walk (QW), which is known as the quantum counterpart model of RW. For that reason, analytical methods for QW are applied to CRW in recent days \cite{afsst,kks1}. In the past two decades, QW has been intensively studied from several points of view \cite{kitagawa,konno1,SMTss} because of its useful applications \cite{kempe, kendon, manouche, portugal}. In the previous study, QW on the path and the cycle graphs are analyzed and the spectra and the eigen functions of the time evolution operators are obtained \cite{araisan, takumikun}.

In this paper, we focus on CRW on the path graph. At first, we consider the spectrum of a Jacobi matrix given by the probability transition matrix (or the time evolution operator) of CRW. This Jacobi matrix is something like the probability transition matrix of RW on the path. 

Next, we discuss a spectral decomposition of the time evolution operator of CRWs with isospectral coins. We should note that under assumption \ref{assump}, all the coins are isospectral but need not be the same because it allows different eigenvectors. We can construct most of the eigenvalues and eigenvectors of the time evolution operator of the CRWs from that of the Jacobi matrix. Although this method is applicable to QWs \cite{araisan,takumikun}, it is not straightforward to apply it to CRW. On the other hand, we always have the eigenvalue $1$ and the corresponding eigenvector of the time evolution operator of CRW, which do not come from that of the Jacobi matrix. This eigenvector corresponds to the birth and death chains and the correspondence means that CRW on the path is supported by RW on the path(birth and death chains). 

The rest of this paper is organized as follows. In Section \ref{defini}, we present the definition of CRW on the path graph and the limit distribution of it. Section \ref{jacobisec} is devoted to the definition and analyzing of the Jacobi matrix given by the time evolution operator of CRW. The main results of this paper, a spectral analysis of CRW are stated in Section \ref{specanal}. 
\section{Definition of CRW on the Path}\label{defini}
In this paper, we consider CRW on the path $P_{n+1}$ with the vertex set $V_{n+1}=\{0,1,\dots,n\}$ and the Edge set $E_{n+1}=\{(x,x+1):x=0,1,\dots,n-1\}$. Next we put the state space ${\mathcal S}_{n+1}={\rm Span}\{|x\rangle\otimes|L\rangle,|x\rangle \otimes|R\rangle:x\in V_{n+1}\}$, where each elements of the tensor product $\{|x\rangle:x\in V_{n+1}\}$ and $\{|L\rangle=\ ^T[1\ 0],|R\rangle=\ ^T[0\ 1]\}$ are orthonormal bases, with $^TA$ represents the transpose of a matrix $A$. Then we define the time evolution operator $U=SC$ on ${\mathcal S}_{n+1}$ with the coin operator $C$ and the shift operator $S$ as follows:
\begin{align*}
C&=\sum_{x=0}^{n}|x\rangle\langle x|\otimes C_x,\\
S|x,J\rangle &=\begin{cases}
|x+1,L\rangle\quad {\rm if}\ x\neq n,\ J=R,\\
|n,R\rangle\quad {\rm if}\ x=n,\ J=R,\\
|x-1,R\rangle\quad {\rm if}\ x\neq 0,\ J=L,\\
|0,L\rangle\quad {\rm if}\ x= 0,\ J=L,
\end{cases}
\end{align*}
where $C_x(x=0,\dots,n)$ are $2\times2$ probability transition matrix and $|x,J\rangle=|x\rangle\otimes|J\rangle$. We should note that the operator $U=SC$ is multiplied from the left. $C_x$ is defined by
\begin{align*}
C_x&=p_{x,L}|L\rangle\langle L|+(1-p_{x,L})|R\rangle\langle L|+p_{x,R}|L\rangle\langle R|+(1-p_{x,R})|R\rangle\langle R|\\
&=\begin{bmatrix}
p_{x,L}&p_{x,R}\\1-p_{x,L}&1-p_{x,R}
\end{bmatrix},
\end{align*}
where $0\leq p_{x,J}\leq 1\ (J=L,R)$. $p_{x,L}$ denotes the probability that the random walker, having moved from $x-1$ to $x$, moves from $x$ to $x-1$. On the other hand, $p_{x,R}$ denotes the probability that the random walker, having moved from $x+1$ to $x$, moves from $x$ to $x-1$. To prevent this model from becoming equivalent to a random walk or becoming a trivial model, we set $0<|p_{x,L}-p_{x,R}|<1$. Next we put $|\Psi_t\rangle\in {\mathcal S}_{n+1}$ be the states of the whole system at time $t\in{\mathbb Z}_{\geq 0}$. In particular, the initial state is described as $|\varphi\rangle=|\Psi_0\rangle$. Let $X_t\in V_{n+1}$ be the position of the random walker at time $t$. Then the probability ${\mathbb P}_{|\varphi\rangle}(X_t=x)$ that the walker with the initial state $|\varphi\rangle$ is located at the vertex $x$ at time $t$ is given by
\begin{align*}
{\mathbb P}_{|\varphi\rangle}(X_t=x)=\Bigl\{\langle x|\otimes\bigl(\langle L|+\langle R|\bigr)\Bigr\} |\Psi_t\rangle.
\end{align*}
We define the limiting distribution $p_{\infty,n}$ of CRW on the path $P_{n+1}$ is described as follows: 
\begin{align*}
p_{\infty,n}(x)=\lim_{t\to\infty}{\mathbb P}_{|\varphi\rangle}(X_t=x),
\end{align*}
where $x=0,\dots,n$.

\section{Jacobi matrix}\label{jacobisec}
In this section, we consider the Jacobi matrix, which is helpful for analyzing CRW. Let $\nu_{1,x},\nu_{2,x}$ and $|w_{1,x}\rangle,\langle q_{1,x}|, |w_{2,x}\rangle,\langle q_{2,x}|$ be the eigenvalues and the corresponding right and left eigenvectors of $C_x$. By direct calculation, we have
\begin{align}
&\nu_{1,x}=1,\quad |w_{1,x}\rangle=\frac{1}{1-(p_{x,L}-p_{x,R})}\begin{bmatrix}p_{x,R}\\1-p_{x,L}
\end{bmatrix},\ \langle q_{1,x}|=\begin{bmatrix}1&1\end{bmatrix},\notag\\
&\nu_{2,x}=p_{x,L}-p_{x,R},\quad |w_{2,x}\rangle=\frac{1}{\sqrt{2}}\begin{bmatrix}1\\-1 \end{bmatrix},\notag\\&\langle q_{2,x}|=\frac{\sqrt{2}}{1-(p_{x,L}-p_{x,R})}
\begin{bmatrix}
1-p_{x,L}&-p_{x,R}
\end{bmatrix}.\label{eigensofc}
\end{align}
Then we consider the spectral decomposition of each matrix $C_x$ as follows:
\begin{align*}
C_x&=1\cdot |w_{1,x}\rangle\langle q_{1,x}|+\nu_{2,x}|w_{2,x}\rangle\langle q_{2,x}|
\\&=(\nu_{2,x}-1)|w_{2,x}\rangle\langle q_{2,x}|+I_2\\
&=(1-\nu_{2,x})|w_{2,x}\rangle\langle -q_{2,x}|+I_2,
\end{align*}
where $I_k$ is the $k\times k$ identity matrix. Here we use the relation $I_2=|w_{1,x}\rangle\langle q_{1,x}|+|w_{2,x}\rangle\langle q_{2,x}|
$ coming from the diagonalizability of $C_x$. This shows that we can represent $C_x$ without $|w_{1,x}\rangle$. 

We define the $(n+1)\times (n+1)$ Jacobi matrix $B$ as follows:
\begin{align*}
\displaystyle{(B)_{i,j}=\begin{cases}
\frac{q_i(L)}{\sqrt{2}}\quad &{\rm if }\ j=i-1,\\
\frac{-q_i(R)}{\sqrt{2}}\quad &{\rm if }\ j=i+1,\\
\frac{-q_0(L)}{\sqrt{2}}\quad &{\rm if }\ i=j=0,\\
\frac{q_n(R)}{\sqrt{2}}\quad &{\rm if }\ i=j=n,\\
0\quad &{\rm otherwise},
\end{cases}}
\end{align*}
where $\langle q_{2,x}|=\begin{bmatrix}q_x(L)&q_x(R) \end{bmatrix}$. We should note that $q_x(L)-q_x(R)=\sqrt2$. In this setting, the corresponding Jacobi matrix is the following:
\begin{tiny}
\begin{align*}
B
=\frac{1}{\sqrt{2}}\begin{bmatrix}
-q_{0}(L)&-q_{0}(R)&&&&&&\\
q_{1}(L)&0&-q_{1}(R)\\
&\ddots&\ddots&\ddots\\
&& q_{x-1}(L)&0&-q_{x-1}(R)&&& \\
&&&q_{x}(L)&0&-q_{x}(R)&& \\
&&&&q_{x+1}(L)&0&-q_{x+1}(R)& \\
&&&&\ddots&\ddots&\ddots\\
&&&&&q_{n-1}(L)&0&-q_{n-1}(R)\\
&&&&&&q_{n}(L)&q_{n}(R)
\end{bmatrix}.
\end{align*}
\end{tiny}
This Jacobi matrix represents the following inner products:
\begin{align*}
\bigl(\langle q_{2,x-1}|\otimes\langle x-1|\bigr)U\bigl(|x\rangle\otimes|w_{2,x}\rangle \bigr)&=\nu_{2,x}(B)_{x-1,x},\\
\bigl(\langle q_{2,x+1}|\otimes\langle x+1|\bigr)U\bigl(|x\rangle\otimes|w_{2,x}\rangle \bigr)&=\nu_{2,x}(B)_{x+1,x},\\
\bigl(\langle q_{2,0}|\otimes\langle 0|\bigr)U\bigl(|0\rangle\otimes|w_{2,0}\rangle \bigr)&=\nu_{2,0}(B)_{0,0},\\
\bigl(\langle q_{2,n}|\otimes\langle n|\bigr)U\bigl(|n\rangle\otimes|w_{2,n}\rangle \bigr)&=\nu_{2,n}(B)_{n,n}.
\end{align*}
This properties help our spectral analysis in Sect.\ref{specanal}. 

$B$ is something like the probability transition matrix of a random walk with a right-hand side operation, as the sum of numbers in the each row (not column) is $1$, except for the $0$th and the $n$th rows. In the $0$th and the $n$th rows, where $-q_0(L)<0$ and $q_n(R)<0$, it can be interpreted as random walkers being absorbed at the end points of the path. 

Next we set an unit vector $\pi=\begin{bmatrix}\pi(0)&\dots&\pi(n) \end{bmatrix}$ such that
\begin{align*}
\pi(0)&=C_{\pi},\\
\pi(x)&=\frac{\prod_{y=1}^{x}q_y(L)}{\prod_{y=0}^{x-1}(-q_y(R))}C_{\pi},
\end{align*}
where $C_{\pi}=\bigl(1+\sum_{x=1}^{n} \frac{\prod_{y=1}^{x}q_y(L)}{\prod_{y=0}^{x-1}(-q_y(R))}\bigr)^{-1}$. We put another Jacobi matrix $J$ as follows:
\begin{scriptsize}
\begin{align*}
&J=D_{\pi^{1/2}}^{-1}BD_{\pi^{1/2}}\\
&=\frac{1}{\sqrt{2}}\begin{bmatrix}
-q_0(L)&\sqrt{-q_0(R)q_1(L)}\\
\sqrt{-q_0(R)q_1(L)}&0&\ddots\\
&\ddots&\ddots&\ddots\\
&&\ddots&0&\sqrt{-q_{n-1}(R)q_{n}(L)}\\
&&&\sqrt{-q_{n-1}(R)q_{n}(L)}&q_n(R)
\end{bmatrix},
\end{align*}
\end{scriptsize}where $D_{\pi^{1/2}}={\rm diag}\Bigl(\sqrt{\pi(0)},\dots,\sqrt{\pi(n)}\Bigr)$ and $A^{-1}$ is the inverse matrix of matrix $A$. Then we have the following proposition.
\begin{proposition}\label{bn1eig}
$B$ has $n+1$ simple real eigenvalues.
\end{proposition}
To prove this proposition, we consider the following lemma.
\begin{lemma}\label{bjiso}
$B$ and $J$ are isospectral.
\end{lemma}
{\it proof of Lemma \ref{bjiso}.} Let $\lambda_m(m=1,\dots,n+1)$ be an eigenvalue of $B$ and ${\bf v}_{m}$ be the corresponding eigenvector. Then $\lambda$ and $v_{m}$ satisfy $B{\bf v}_{m}=\lambda_m {\bf v}_{m}$. Noting that $J=D_{\pi^{1/2}}^{-1}BD_{\pi^{1/2}}$, we have
\begin{align*}
D^{-1}_{\pi^{1/2}}B{\bf v}_{m}&=D^{-1}_{\pi^{1/2}}\lambda_m {\bf v}_{m},\\
D^{-1}_{\pi^{1/2}}BD_{\pi^{1/2}}D^{-1}_{\pi^{1/2}}{\bf v}_{m}&=\lambda_m D^{-1}_{\pi^{1/2}} {\bf v}_{m},\\
JD^{-1}_{\pi^{1/2}}{\bf v}_{m}&=\lambda_m D^{-1}_{\pi^{1/2}} {\bf v}_{m}.
\end{align*}
Thus $\lambda_m$ is an eigenvalue of $J$ and $D^{-1}_{\pi^{1/2}} {\bf v}_{m}$ is the corresponding eigenvector. Since this correspondence holds for any eigenvalue of $B$, ${\rm Spec}(B) \subset {\rm Spec}(J)$. The similar argument holds for the inverse as well.\hfill $\square$

By general argument about the Jacobi matrix, all of the eigenvalues of $J$ are simple. Combining this with Lemma \ref{bjiso}, we obtain Proposition \ref{bn1eig}. 
\section{A spectral analysis of CRW on the Path}\label{specanal}
In this section, we give a framework of spectral analysis for CRW on $P_{n+1}$. In order to do so, we restrict the coin operator as follows:
\begin{assumption}\label{assump}
We assume that the coin operator consists of isospectral probability transition matrices, i.e., we use
\begin{align*}
C=\sum_{x=0}^{n}|x\rangle\langle x|\otimes \Bigl\{(1-\nu_{2})|w_{2,x}\rangle\langle -q_{2,x}|+I_2\Bigr\},
\end{align*}
as the coin operator, where $\nu_2=p_{x,L}-p_{x,R}$ for all $x=0,\dots,n$.
\end{assumption}
\subsection{The spectrum of $U$ from $B$}\label{subset1}
This subsection deals with the eigenvalues and the corresponding eigenvectors of $U$ related to those of $B$.

Let $\lambda_1>\dots>\lambda_{n+1}$ be the eigenvalues and ${\bf v}_{m}=\ ^T\begin{bmatrix}v_m(0)&\dots&v_m(n) \end{bmatrix}$ $(m=1,\dots,n+1)$ be the corresponding normalized eigenvectors of $B$. For each $\lambda_m$ and ${\bf v}_m$, we define two vectors
\begin{align*}
{\bf a}_m &=\sum_{x=0}^{n}v_m(x)\bigl(|x\rangle\otimes|w_{2}\rangle\bigr)\\
&=\frac{1}{\sqrt{2}}\sum_{x=0}^{n}\Bigl(v_m(x)|x\rangle\otimes\bigl(|L\rangle-|R\rangle\bigr)\Bigr),\\
{\bf b}_m&=S{\bf a}_m\\
&=\frac{1}{\sqrt{2}}\sum_{x=1}^{n-1}|x\rangle\otimes\Bigl(v_m(x+1)|R\rangle-v_m(x-1)|L\rangle\Bigr)\\
&+\frac{|0\rangle}{\sqrt{2}}\otimes\Bigl(v_m(1)|R\rangle+v_m(0)|L\rangle\Bigr)+\frac{|n\rangle}{\sqrt{2}}\otimes\Bigl(-v_m(n)|R\rangle-v_m(n-1)|L\rangle\Bigr).
\end{align*}
By using $S^2=I_{n+1}\otimes I_2$, it is easy to see that $C{\bf a}_m=\nu_2{\bf a}_m$ and then $U{\bf a}_m=\nu_2{\bf b}_m$. Also we have $C{\bf b}_m=(1-\nu_2)\lambda_m{\bf a}_m+{\bf b}_m$ and $U{\bf b}_m={\bf a}_m+(1-\nu_2)\lambda_m{\bf b}_m$. So we have the following relationship:
\begin{align}
U\begin{bmatrix}
{\bf a}_m\\{\bf b}_m
\end{bmatrix}=
\begin{bmatrix}
0&\nu_2\\
1&(1-\nu_2)\lambda_m
\end{bmatrix}
\begin{bmatrix}
{\bf a}_m\\{\bf b}_m
\end{bmatrix}.\label{linearab}
\end{align}
We also obtain $\|{\bf a}_m\|=\|{\bf b}_m\|=1$. Let $({\bf a}_m,{\bf b}_m)$ be the inner product of vectors ${\bf a}_m$ and ${\bf b}_m$. When $|({\bf a}_m,{\bf b}_m)|<1$, ${\bf a}_m$ and ${\bf b}_m$ are linearly independent. About the linearly independence, we have the following lemma.
\begin{lemma}\label{linedep}
The vectors ${\bf a}_{m}$ and ${\bf b}_{m}=S{\bf a}_{m}$ which are related to $\lambda_m$ are linearly dependent if and only if $\lambda_{n+1}=-1$.
\end{lemma}
{\it proof of Lemma \ref{linedep}.} By direct calculation, we have 
\begin{align*}
({\bf a}_m,{\bf b}_m)&=\frac{1}{2}\sum_{x=1}^{n-1}\bigl(-v_m(x-1)v_m(x)-v_m(x)v_m(x+1) \bigr) \\
&+\frac{1}{2}\bigl(|v_m(0)|^2-v_m(0)v_m(1)\bigr)+\frac{1}{2}\bigl(-v_m(n-1)v_m(n)+|v_m(n)|^2\bigr)\\
&=-\sum_{x=1}^{n}\Bigl(v_m(x-1)v_m(x)\Bigr)+\frac{1}{2}\bigl(|v_m(0)|^2+|v_m(n)|^2\bigr).
\end{align*}
By triangle inequality and Cauchy-Schwarz inequality, we get
\begin{align}
|({\bf a}_m,{\bf b}_m)|&\leq \sum_{x=1}^{n}|v_m(x-1)v_m(x)|+\frac{1}{2}\bigl(|v_m(0)|^2+|v_m(n)|^2\bigr)\notag\\
&\leq \sqrt{\Biggl(\sum_{x=0}^{n-1}|v_m(x)|^2\Biggr)\Biggl(\sum_{x=1}^{n}|v_m(x)|^2\Biggr)}+\frac{1}{2}\bigl(|v_m(0)|^2+|v_m(n)|^2\bigr)\notag\\
&=\sqrt{\bigl(1-|v_m(n)|^2\bigr)\bigl(1-|v_m(0)|^2\bigr)}+\frac{1}{2}\bigl(|v_m(0)|^2+|v_m(n)|^2\bigr).\label{inequalities}
\end{align}
Since the rightmost side of inequality \eqref{inequalities} takes the maximum value $1$ when $|v_m(0)|=|v_m(n)|$, we have
\begin{align}
|({\bf a}_m,{\bf b}_m)|\leq 1.\label{fulineq}
\end{align}
Although inequality \eqref{fulineq} is trivial, the conditions for equality in inequalities \eqref{inequalities} and \eqref{fulineq} are important. By the equality condition of triangle inequality in \eqref{inequalities}, we obtain 
\begin{align}
v_m(x-1)v_m(x)\leq 0\quad(x=1,\dots,n).\label{eqcon1}
\end{align}
The equality condition of Cauchy-Schwarz inequality gives
\begin{align}
v_m(x-1)v_m(y)=v_m(x)v_m(y-1),\label{eqcon2}
\end{align}
where $x\neq y$ and $x,y=1,\dots, n$. Combining $|v_m(0)|=|v_m(n)|$ with conditions \eqref{eqcon1} and \eqref{eqcon2}, we have
\begin{align*}
v_m(x)=-v_m(x-1)\quad (x=1,\dots,n).
\end{align*}
Owing to $\|{\bf v}_m\|=1$, we obtain
\begin{align}
v_m(0)=\pm\frac{1}{\sqrt{n+1}},v_m(1)=\mp \frac{1}{\sqrt{n+1}},\dots,\label{fukugoudoujun}
\end{align}
where the double-sign corresponds. Because both of the case in \eqref{fukugoudoujun} satisfy $B{\bf v}_m=-{\bf v}_m$, this is an eigenvector of $B$ corresponding to $\lambda_m=-1$. Due to the simplicity of the eigenvalues of $B$, we have the uniqueness of $\lambda_m=-1$ for the linear dependence between ${\bf a}_m$ and ${\bf b}_m$.\hfill $\square$

As we see in sections \ref{nu2seisec} and \ref{nu2fusec}, $|\lambda_m|\leq 1$. Thus, the minimum eigenvalue $\lambda_{n+1}$ of $B$ denotes $\lambda_{n+1}=-1$.

For the case with $\lambda_{n+1}=-1$, we obtain an eigenvalue 
$\mu_{n+1}=\nu_2$ of $U$ and the corresponding eigenvector ${\bf a}_{n+1}={\bf b}_{n+1}$ by direct calculation.

For the other cases with $\lambda_m\neq -1$, we see from Eq. \eqref{linearab} that $U$ is a linear operator acting on the linear space, ${\rm Span}({\bf a}_m, {\bf b}_m)$. In order to obtain the eigenvalues and eigenvectors, we take a vector $\alpha{\bf a}_m+\beta{\bf b}_m\in {\rm Span}({\bf a}_m, {\bf b}_m)$. The eigen equation for $U$ is given by $U(\alpha{\bf a}_m+\beta{\bf b}_m)=\mu(\alpha{\bf a}_m+\beta{\bf b}_m)$. From Eq. \eqref{linearab}, this is equivalent to
\begin{align*}
    \begin{bmatrix}
0&1\\ \nu_2&(1-\nu_2)\lambda_m
    \end{bmatrix}
    \begin{bmatrix}
\alpha\\ \beta
\end{bmatrix}
=\mu\begin{bmatrix}
\alpha\\ \beta
\end{bmatrix}.
\end{align*}
Therefore we can obtain two eigenvalues $\mu_{\pm m}$ of $U$ which are related to the eigenvalue $\lambda_m$ of $B$ as solutions of the following quadratic equation:
\begin{align}
\mu_{\pm m}^2-(1-\nu_2)\lambda_m\mu_{\pm m}-\nu_2=0.\label{eigenequ}
\end{align}
Also we have the corresponding eigenvectors ${\bf u}_{\pm m}={\bf a}_m+\mu_{\pm m}{\bf b}_m$ by setting $\alpha=1,\beta=\mu_{\pm m}$. By direct calculation, we obtain
\begin{align}
\mu_{\pm m}=\frac{(1-\nu_2)\lambda_{m}\pm\sqrt{(1-\nu_2)^2\lambda_{m}^2+4\nu_2}}{2}. \label{mupmatai}
\end{align}
In \cite{afsst}, for the spatially homogeneous coin case, these eigenvalues $\mu_{\pm m}$ were obtained.

When $\lambda_m=-1$, we have $\mu_{+m}=\nu_2$ and $\mu_{-m}=-1$ from Eq.\eqref{mupmatai}. However, in the case $\lambda_m=-1$, we know that ${\bf a}_m$ and ${\bf b}_m$ are linearly dependent, and we obtain only one eigenvalue $\mu_m=\nu_2$ and the corresponding eigenvector ${\bf a}_m$ of $U$. Connecting this with the simplicity of $\mu_{\pm m}$ with respect to $\lambda_m$ as we see later in sections \ref{nu2seisec} and \ref{nu2fusec}, we have $\mu_{\pm m}\neq -1$. 

The relationship between solution and coefficients of Eq. \eqref{eigenequ} gives
\begin{align}
\begin{cases}
\mu_{-m}+\mu_{+m}=(1-\nu_2)\lambda_m,\\
\mu_{-m}\mu_{+m}=-\nu_2.
\end{cases}\label{kaitokeisu}
\end{align}
Then we deal with the following two cases: \ref{nu2seisec} $\nu_2>0$ and \ref{nu2fusec} $\nu_2<0$. 
\newpage
\subsubsection{$\nu_2>0$ case}\label{nu2seisec}

By Eq. \eqref{kaitokeisu}, we have the following figure, where $\mu_{\pm m}$ are the $x$-coordinates of the intersections of the straight line and the inverse proportional curve.

\begin{figure}[htbp]
\begin{center}
\begin{tikzpicture}
\coordinate[label=below left:O](O)at(0,0);
\coordinate(XS)at(-4,0);
\coordinate(XL)at(4,0);
\coordinate(YS)at(0,-4);
\coordinate(YL)at(0,4);
\draw[name path=xax][semithick,->,>=stealth](XS)--(XL)node[right]{$x$};
\draw[name path=yax][semithick,->,>=stealth](YS)--(YL)node[above]{$y$};
\draw[name path=sekis][thick,samples=100,domain=-4:-0.25]plot(\x,-1/\x);
\draw[name path=seki][thick,samples=100,domain=0.25:4]plot(\x,-1/\x)node[below]{$\displaystyle y=\frac{-\nu_2}{x}$};
\draw[name path=wa][thick,domain=-2.5:4]plot(\x,-\x+1.5)node[below]{$y=-x+(1-\nu_2)\lambda_m$};
\draw[name path=waw][dashed,domain=-0.5:4]plot(\x,-\x+3.5)node[above]{};
\draw[name path=wawa][dashed,domain=-4:0.5]plot(\x,-\x-3.5)node[above]{};
\path[name intersections={of= seki and wa, by={beta}}];
\path[name intersections={of= sekis and wa, by={alpha}}];
\path[name intersections={of= sekis and waw, by={alpp}}];
\path[name intersections={of= sekis and wawa, by={alm}}];
\path[name intersections={of= seki and waw, by={bem}}];
\path[name intersections={of= seki and wawa, by={bep}}];
\coordinate[label=right:$(1-\nu_2)\lambda_m$](nulam)at(0,1.5);
\coordinate[label=right:$1-\nu_2$](nup)at(0,3.5);
\coordinate[label=left:$-(1-\nu_2)$](num)at(0,-3.5);
\draw[name path=xalpha][dashed]($(XS)!(alpha)!(XL)$)node[above left]{$\mu_{-m}$}--(alpha);
\draw[name path=xbeta][dashed]($(YS)!(alpha)!(YL)$)node[right]{$\mu_{+m}$}--(alpha);
\draw[name path=yalpha][dashed]($(XS)!(beta)!(XL)$)node[above]{$\mu_{+m}$}--(beta);
\draw[name path=ybeta][dashed]($(YS)!(beta)!(YL)$)node[left]{$\mu_{-m}$}--(beta);
\draw[name path=alphap][dashed]($(XS)!(alpp)!(XL)$)node[above]{\begin{small}$-\nu_2$\end{small}}--(alpp);
\draw[name path=alpham][dashed]($(XS)!(alm)!(XL)$)node[below]{$-1$}--(alm);
\draw[name path=betam][dashed]($(XS)!(bem)!(XL)$)node[above]{$1$}--(bem);
\draw[name path=betap][dashed]($(XS)!(bep)!(XL)$)node[above]{$\nu_2$}--(bep);
\path[name intersections={of= xax and xalpha, by={tal}}];
\path[name intersections={of= yax and xbeta, by={tbe}}];
\path[name intersections={of= xax and yalpha, by={ttal}}];
\path[name intersections={of= yax and ybeta, by={ttbe}}];
\path[name intersections={of= xax and alphap, by={alphamax}}];
\path[name intersections={of= xax and alpham, by={alphamin}}];
\path[name intersections={of= xax and betap, by={betamax}}];
\path[name intersections={of= xax and betam, by={betamin}}];
\draw[red][thick](alphamin)--(alphamax);
\draw[blue][thick](betamin)--(betamax);
\foreach\P in{alpha,beta,nulam,nup,num,alpp,alm,bem,bep,tal,tbe,ttal,ttbe,betamax}\fill[black](\P)circle(0.05); 
\foreach\P in{alphamin,alphamax,betamin}\filldraw[fill=white,draw=black](\P)circle(0.06);
\end{tikzpicture}
\caption{Mapping from $\lambda_m$ to $\mu_{\pm m}\ (\nu_2>0)$}\label{nu2sei}
\end{center}
\end{figure}
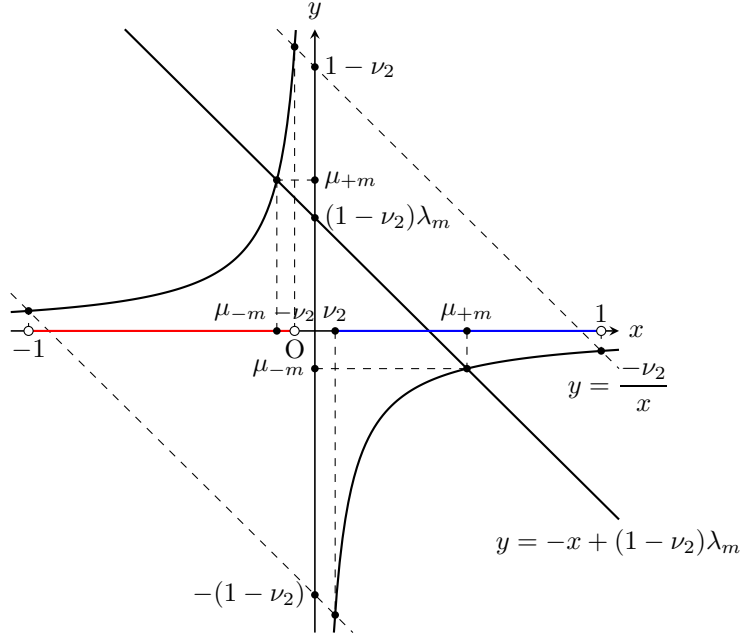

Since $U$ is the probability transition matrix, $|\mu_{\pm m}|\leq 1$, and we also observe from Fig. \ref{nu2sei} that $\mu_{\pm m}$ increases monotonically with respect to $\lambda_m$, which determines the $y-$intercept $(1-\nu_2)\lambda_m$ of the line $y=-x+(1-\nu_2)\lambda_m$. Thus we have $|\lambda_m|\leq 1$ and every element in ${\rm Spec}(U)$ which we obtain from $B$ is simple. Putting it all together, we have $2n+1$ eigenvalues $\mu_{\pm m}$ and corresponding eigenvectors ${\bf u}_{\pm m}$ of $U$, which are related to $n+1$ eigenvalues $\lambda_m$ and corresponding eigenvectors ${\bf v}_m$ of $B$.

In addition, we obtain $\lambda_m\neq1$ in Sec. \ref{anotherspe}.
\newpage
\noindent
\subsubsection{$\nu_2<0$ case}\label{nu2fusec}

By Eq. \eqref{kaitokeisu}, we have the following figure.

\begin{figure}[htbp]
\begin{center}
\begin{tikzpicture}
\coordinate[label=below left:O](O)at(0,0);
\coordinate(XS)at(-4,0);
\coordinate(XL)at(4,0);
\coordinate(YS)at(0,-4);
\coordinate(YL)at(0,4);
\draw[name path=xax][semithick,->,>=stealth](XS)--(XL)node[right]{$x$};
\draw[name path=yax][semithick,->,>=stealth](YS)--(YL)node[above]{$y$};
\draw[name path=sekis][thick,samples=100,domain=0.25:4]plot(\x,1/\x);
\draw[name path=seki][thick,samples=100,domain=-4:-0.25]plot(\x,1/\x)node[below]{$\displaystyle y=\frac{-\nu_2}{x}$};
\draw[name path=wa][thick,domain=-0.8:4]plot(\x,-\x+2.6)node[below]{$y=-x+(1-\nu_2)\lambda_m$};
\draw[name path=waa][thick,domain=-1.4:-1.1]plot(\x,-\x+2.6)node[]{};
\draw[name path=waw][dashed,domain=-0.5:4]plot(\x,-\x+3.5)node[above]{};
\draw[name path=wawa][dashed,domain=-4:0.5]plot(\x,-\x-3.5)node[above]{};
\draw[name path=walim][dashed,domain=-4:1.3]plot(\x,-\x-2)node[above]{};
\draw[name path=walimowa][dashed,domain=1.7:2]plot(\x,-\x-2)node[above]{};
\draw[name path=walip][dashed,domain=-0.4:2]plot(\x,-\x+2)node[]{};
\draw[name path=walipp][dashed,domain=2.4:3.5]plot(\x,-\x+2)node[above]{};
\draw[name path=walippp][dashed,domain=3.9:4]plot(\x,-\x+2)node[]{};
\draw[name path=walippa][dashed,domain=-2:-0.8]plot(\x,-\x+2)node[above]{};
\path[name intersections={of= sekis and wa, by={alpha,beta}}];
\path[name intersections={of= sekis and waw, by={alpm,alpp}}];
\path[name intersections={of= yax and waw, by={seppe}}];
\path[name intersections={of= sekis and wawa, by={alm,almp}}];
\path[name intersections={of= yax and wawa, by={seppem}}];
\path[name intersections={of= seki and walim, by={doublem}}];
\path[name intersections={of= sekis and walip, by={doublep}}];
\path[name intersections={of= seki and wawa, by={bem,bep}}];
\coordinate[label=left:$(1-\nu_2)\lambda_m$](nulam)at(0,2.6);
\coordinate[label=left:$1-\nu_2$](nup)at(0,3.5);
\coordinate[label=right:$-(1-\nu_2)$](num)at(0,-3.5);
\draw[name path=xalpha][dashed]($(XS)!(alpha)!(XL)$)node[above right=-0.06cm]{\begin{small}$\mu_{-m}$\end{small}}--(alpha);
\draw[name path=yalpha][dashed]($(XS)!(beta)!(XL)$)node[below]{$\mu_{+m}$}--(beta);
\draw[name path=alphapm][dashed]($(XS)!(alpm)!(XL)$)node[below]{\begin{small}$-\nu_2$\end{small}}--(alpm);
\draw[name path=alphap][dashed]($(XS)!(alpp)!(XL)$)node[below]{$1$}--(alpp);
\draw[name path=alpham][dashed]($(XS)!(alm)!(XL)$)node[above]{}--(alm);
\draw[name path=betam][dashed]($(XS)!(bem)!(XL)$)node[above]{$-1$}--(bem);
\draw[name path=betap][dashed]($(XS)!(bep)!(XL)$)node[above]{$\nu_2$}--(bep);
\draw[name path=doubmx][dashed]($(XS)!(doublem)!(XL)$)node[above left=-0.1cm]{\begin{small}$-\sqrt{-\nu_2}$\end{small}}--(doublem);
\draw[name path=doubpx][dashed]($(XS)!(doublep)!(XL)$)node[below]{\begin{small}$\sqrt{-\nu_2}$\end{small}}--(doublep);
\path[name intersections={of= xax and xalpha, by={tal}}];
\path[name intersections={of= xax and yalpha, by={tbe}}];
\path[name intersections={of= yax and yalpha, by={tttbe}}];
\path[name intersections={of= xax and alphap, by={alphamax}}];
\path[name intersections={of= xax and alphapm, by={alphamin}}];
\path[name intersections={of= xax and betap, by={betamax}}];
\path[name intersections={of= xax and betam, by={betamin}}];
\path[name intersections={of= xax and doubmx, by={doumbor}}];
\path[name intersections={of= xax and doubpx, by={doupbor}}];
\draw[blue][thick](doupbor)--(alphamax);
\draw[red][thick](alphamin)--(doupbor);
\draw[blue][thick](doumbor)--(betamax);
\draw[red][thick](betamin)--(doumbor);
\foreach\P in{alpha,beta,alpm,alpp,alm,almp,betamax,doublem,doublep,doumbor,doupbor,tal,tbe,nulam,tbe,tttbe,alpm,alpp,bem,bep,seppe,seppem}\fill[black](\P)circle(0.05); 
\foreach\P in{alphamin,alphamax,betamin}\filldraw[fill=white,draw=black](\P)circle(0.06);
\end{tikzpicture}
\caption{Mapping from $\lambda_m$ to $\mu_{\pm m}\ (\nu_2<0)$}\label{nu2fu}
\end{center}
\end{figure}
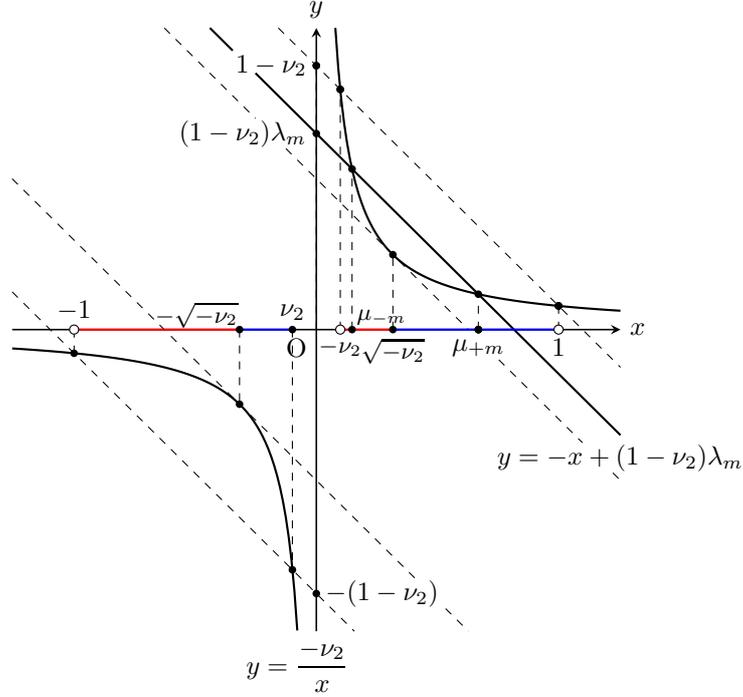
As we see in Sec. \ref{nu2seisec}, each of $\mu_{+ m}$ and $\mu_{-m}$ is simple with respect to $\lambda_m$.

In this case, as we observe from Fig. \ref{nu2fu}, some eigenvalues $\mu_{\pm m}$ of $U$ could be the double root $\mu_m=\pm\sqrt{-\nu_2}$ or the imaginary root of the eigen equation. Thus we consider the following assumption.
\begin{assumption}\label{assump2}
    We assume that all of the eigenvalues $\lambda_m$ of $B$ satisfy $\frac{\sqrt{-4\nu_2}}{1-\nu_2}<|\lambda_m|\leq 1$.
\end{assumption}
Under this assumption, we have $2n+1$ eigenvalues $\mu_{\pm m}$ and corresponding eigenvectors ${\bf u}_{\pm m}$ of $U$, which are related to $n+1$ eigenvalues $\lambda_m$ and corresponding eigenvectors ${\bf v}_m$ of $B$. 

Furthermore, we obtain $\lambda_m\neq1$ in Sec. \ref{anotherspe}.
\newpage
\subsection{Another spectrum of $U$}\label{anotherspe}

In this subsection, we consider the eigenvalue and the corresponding eigenvector of $U$ which is not related to $B$. We should note that we already have $2n+1$ eigenvalues and eigenvectors from $B$, and $U$ is a $2(n+1)\times2(n+1)$ probability transition matrix. Thus, we still need one eigenvalue and the corresponding eigenvector of $U$. 

We consider the following vector:
\begin{align*}
{\bf u}_0&=\sum_{x=0}^{n}v_0(x)|x\rangle\otimes\frac{1}{1-\nu_2}\begin{bmatrix}
p_{x,R}\\1-p_{x,L}
\end{bmatrix},\\
&=\sum_{x=0}^{n}v_0(x)|x\rangle\otimes|w_{1,x}\rangle,
\end{align*}
where $\sum_{x=0}^{n}v_0(x)=1$. Noting that all of the eigenvectors of $U$ obtained from $B$ satisfy
\begin{align}
\sum_{x=0}^{n}\Bigl\{\langle x|\otimes\bigl(\langle L|+\langle R|\bigr)\Bigr\}{\bf u}_{m}=0,\label{abnopro}
\end{align}
where $m\in\{\pm1,\dots,\pm n,n+1\}$, we have ${\bf u}_0 \notin {\rm Span}({\bf u}_{\pm 1},\dots,{\bf u}_{\pm n},{\bf u}_{n+1})$. Eq. \eqref{abnopro} means that we can not obtain any probability distribution from any states in Span$({\bf u}_{\pm 1},\dots,{\bf u}_{\pm n},{\bf u}_{n+1})$. Since Eq. \eqref{eigensofc} shows $C_x|w_{1,x}\rangle=|w_{1,x}\rangle$ for all $x=0,\dots,n$, ${\bf u}_0$ is the eigenvector of the coin operator $C$, corresponding to the eigenvalue $\nu_1=1$. Therefore when ${\bf u}_0$ also satisfies $S{\bf u}_0={\bf u}_0$, ${\bf u}_0$ is the eigenvector of the time evolution operator $U$ corresponding to the eigenvalue $\mu_0=1$. Thus we have the following conditions:
\begin{align}
v_0(x)(1-p_{x,L})&=v_0(x+1)p_{x+1,R},\notag\\
v_0(x+1)&=\frac{1-p_{x,L}}{p_{x+1,R}}v_{0}(x),\label{detbal}
\end{align}
where $x=0,\dots,n-1$. By Eq. \eqref{detbal}, we obtain
\begin{align*}
v_0(x)=\frac{\prod_{y=0}^{x-1}(1-p_{y,L})}{\prod_{y=1}^{x}p_{y,R}}v_{0}(0),
\end{align*}
where $x=1,\dots,n$. Owing to $\sum_{x=0}^{n}v_0(x)=1$, we have
\begin{align*}
v_0(0)=\Biggl(1+\sum_{x=1}^n\frac{\prod_{y=0}^{x-1}(1-p_{y,L})}{\prod_{y=1}^{x}p_{y,R}}\Biggr)^{-1},
\end{align*}
and $v_0$ can be regarded as the probability distribution. We should note that 
\begin{align}
\sum_{x=0}^{n}\Bigl\{\langle x|\otimes\bigl(\langle L|+\langle R|\bigr)\Bigr\}{\bf u}_0=\sum_{x=0}^{n}v_0(x)=1.\label{uopro}
\end{align}
Let $Q$ be a probability transition matrix described as follows:
\begin{footnotesize}
\begin{align*}
Q=\frac{1}{1-\nu_2}\begin{bmatrix}
p_{0,R}&p_{1,R}\\
1-p_{0,L}&0&p_{2,R}\\
&1-p_{1,L}&0&\ddots\\
&&\ddots&\ddots&\ddots\\
&&&\ddots&0&p_{n-1,R}\\
&&&&1-p_{n-2,L}&0&p_{n,R}\\
&&&&&1-p_{n-1,L}&1-p_{n,L}
\end{bmatrix}.
\end{align*}
\end{footnotesize}
Then, we can regard $v_0$ as the probability distribution of RW on the path whose probability transition matrix is $Q$. We see that $v_0$ satisfies the detailed balance conditions by Eq. \eqref{detbal}. Thus, $v_0$ is the stationary distribution of birth and death chains. Combining this with Eq. \eqref{abnopro}, we consider that CRW on the path is supported by birth and death chains.

Moreover, by the irreducibility of $U$ and the Perron–Frobenius theorem, the largest eigenvalue $\mu_0=1$ of the probability transition matrix $U$ is simple. Therefore all of the eigenvalues $\mu_{\pm m}$ of $U$ obtained from $B$ are not $1$. 

\subsection{Main results}
This subsection is devoted to our main results.

As a consequence of our discussion so far, we obtain the following theorem:

\begin{theorem}\label{mainthe}
If $\nu_2>0$ or $\nu_2<0$ and every $\lambda_m(m=1,\dots,n+1)$ satisfies $\frac{\sqrt{-4\nu_2}}{1-\nu_2}<|\lambda_m|\leq 1$, every element in ${\rm Spec}(U)$ is simple. The eigenvalues $\mu_0,\mu_{\pm m},\mu_{n+1}$ and the corresponding eigenvectors ${\bf u}_0,{\bf u}_{\pm m},{\bf u}_{n+1}$ are following:\\
$(1)\ \mu_0=1$ and 
\begin{align*}
{\bf u}_0=\sum_{x=0}^{n}\frac{\prod_{y=0}^{x-1}(1-p_{y,L})}{\prod_{y=1}^{x}p_{y,R}}v_{0}(0)|x\rangle\otimes\frac{1}{1-\nu_2}\begin{bmatrix}
p_{x,R}\\1-p_{x,L}\end{bmatrix},\\
{\rm where}\ v_{0}(0)=\Biggl(1+\sum_{x=1}^n\frac{\prod_{y=0}^{x-1}(1-p_{y,L})}{\prod_{y=1}^{x}p_{y,R}}\Biggr)^{-1}.
\end{align*}
$(2)$ For $m=1,\dots,n$,
\begin{align*}
\mu_{\pm m}=\frac{(1-\nu_2)\lambda_{m}\pm\sqrt{(1-\nu_2)^2\lambda_{m}^2+4\nu_2}}{2},
\end{align*}
and
\begin{small}
\begin{align*}
{\bf u}_{\pm m}&={\bf a}_m+\mu_{\pm m}{\bf b}_m\\
&=\frac{1}{\sqrt{2}}\sum_{x=1}^{n-1}|x\rangle\otimes\Bigl\{\bigl(v_m(x)-\mu_{\pm m}v_m(x-1)\bigr)|L\rangle+\bigl(\mu_{\pm m}v_m(x+1)-v_m(x)\bigr)|R\rangle\bigr)\Bigr\}\\
&+\frac{|0\rangle}{\sqrt{2}}\otimes\Bigl\{\bigl(1+\mu_{\pm m}\bigr)v_m(0)|L\rangle+\bigl(\mu_{\pm m}v_m(1)-v_m(0)\bigr)|R\rangle\Bigr\}\\
&+\frac{|n\rangle}{\sqrt{2}}\otimes\Bigl\{\bigl(v_m(n)-\mu_{\pm m}v_m(n-1)\bigr)|L\rangle-\bigl(1+\mu_{\pm m}\bigr)v_m(n)|R\rangle\Bigr\}.
\end{align*}
\end{small}
$(3)\ \mu_{n+1}=\nu_2$ and
\begin{align*}
{\bf u}_{n+1}&={\bf a}_{n+1}\\
&=\frac{1}{\sqrt{2}}\sum_{x=0}^{n}\Bigl\{\frac{(-1)^x}{\sqrt{n+1}}|x\rangle\otimes\bigl(|L\rangle-|R\rangle\bigr)\Bigr\}.
\end{align*}
\end{theorem}

Let $\langle {\bf q}_{m}|\ (m=0,\pm 1,\dots,\pm n,n+1)$ be the projection onto ${\bf u}_{m}=|{\bf u}_{m}\rangle$ on ${\mathcal S}_{n+1}$. By Theorem \ref{mainthe} and the linear independence of the eigenvectors, we obtain the following description about the time evolution operator $U$.
\begin{align*}
U&=\sum_{m\in \{0,\pm 1,\dots,\pm n,n+1\}}\mu_m|{\bf u}_m\rangle\langle {\bf q}_m|,\\
U^t&=\sum_{m\in \{0,\pm 1,\dots,\pm n,n+1\}}\mu_m^t|{\bf u}_m\rangle\langle {\bf q}_m|.
\end{align*}
We have the following corollary.
\begin{corollary}\label{corlim}
If $\nu_2>0$ or $\nu_2<0$ and every $\lambda_m(m=1,\dots,n+1)$ satisfies $\frac{\sqrt{-4\nu_2}}{1-\nu_2}<|\lambda_m|\leq 1$, there exists the following limiting distribution $p_{\infty,n}$ of CRW on the path $P_{n+1}$,
\begin{align}
p_{\infty,n}(0)&=\Biggl(1+\sum_{x=1}^n\frac{\prod_{y=0}^{x-1}(1-p_{y,L})}{\prod_{y=1}^{x}p_{y,R}}\Biggr)^{-1},\label{pinft0}\\
p_{\infty,n}(x)&=\frac{\prod_{y=0}^{x-1}(1-p_{y,L})}{\prod_{y=1}^{x}p_{y,R}}p_{\infty,n}(0)\quad(x=1,\dots,n).\label{pinftx}
\end{align}
\end{corollary}
{\it proof of Corollary \ref{corlim}.}  Since the eigenvectors of $U$ are linearly independent, the initial state of CRW $|\varphi\rangle$ can be written as follows:
\begin{align*}
|\varphi\rangle=\sum_{m\in \{0,\pm 1,\dots,\pm n,n+1\}}\alpha_m|{\bf u}_m\rangle,
\end{align*}
where $\alpha_m\in{\mathbb R}$. Because $|\varphi\rangle\in{\mathcal S}_{n+1}$ is the initial state of the probability process, we get $\alpha_0=1$ for any initial state $|\varphi\rangle$, by Eqs. \eqref{abnopro} and \eqref{uopro}. Then we have the following description for the state $|\Psi_t\rangle$ at time $t$,
\begin{align}
|\Psi_t\rangle&=U^t|\varphi\rangle\notag\\
&=1^t\cdot1\cdot|{\bf u}_0\rangle+\sum_{m\in \{\pm 1,\dots,\pm n,n+1\}}\mu_m^t\alpha_m|{\bf u}_m\rangle.\label{psit}
\end{align}
Owing to $|\mu_{\pm m}|<1\ (m=1,\dots,n)$ and $|\mu_{n+1}|<1$, as $t$ increases to infinity, the second term on the right side of the Eq. \eqref{psit} converges to ${\bf 0}$. Therefore, there exists the limiting distribution $p_{\infty,n}$ and that is written by
\begin{align*}
p_{\infty,n}(x)&=\lim_{t\to\infty}{\mathbb P}_{|\varphi\rangle}(X_t=x)\\
&=\Bigl\{\langle x|\otimes\bigl(\langle L|+\langle R|\bigr)\Bigr\}|{\bf u}_0\rangle.
\end{align*}
By direct calculation, we have Eqs. \eqref{pinft0} and \eqref{pinftx}.\hfill $\square$

In \cite{someexpli}, the limiting distribution is given and that is similar to the spatially homogeneous coin case of our model. The difference between the limiting distributions come from the difference of the model settings at the end points of the path.

\section{Summary}
In this paper, we considered CRW on the path graph. We obtained the eigenvalues and the corresponding eigenvectors of the time evolution operator of CRWs with isospectral but need not the same coins. Most of those eigenvalues and eigenvectors were constructed by those of the Jacobi matrix. We also had the rest eigenvalue $1$ and the corresponding eigenvector coming from the reversible distribution, hence the limit distribution was given by the eigenvector. One of the interesting future problems is considering the CRW on more general graphs. 
\section*{Acknowledgement}
This work is partially financed by the Grant-in-Aid for Young Scientists of Japan Society for the Promotion of Science (Grant No. 23K13017).

\end{document}